\documentclass[12pt]{article}
\usepackage{amsmath}
\usepackage{amssymb}
\usepackage{vatola}

\def\q{\quad}
\def\qq{\qquad}
\def\mod#1{\ (\text{\rm mod}\ #1)}
\def\t{\hbox}
\def\qtq#1{\q\t{#1}\q}
\def\f{\frac}
\def\e{\equiv}

\def\b{\binom}

\def\sls#1#2{(\f{#1}{#2})}
 \def\ls#1#2{\big(\f{#1}{#2}\big)}
\def\Ls#1#2{\Big(\f{#1}{#2}\Big)}

\def\qs#1#2{\big(\f{#1}{#2}\big)_4}

\let \pro=\proclaim
\let \endpro=\endproclaim
\begin{document}
\leftline{preprint: December 2, 2013}\par\q\par\q
 \centerline {\bf
Quartic residues and sums involving $\binom{4k}{2k}$}
$$\q$$
\centerline{Zhi-Hong Sun} $$\q$$ \centerline{School of Mathematical
Sciences, Huaiyin Normal University,} \centerline{Huaian, Jiangsu
223001, P.R. China} \centerline{E-mail: zhihongsun@yahoo.com}
\centerline{Homepage: http://www.hytc.edu.cn/xsjl/szh}

\abstract{Let $p$ be an odd prime and let $m\not\equiv 0\pmod p$ be
a rational p-adic integer. In this paper we reveal the connection
between quartic residues and the sum
$\sum_{k=0}^{[p/4]}\binom{4k}{2k}\frac 1{m^k}$, where $[x]$ is the
greatest integer not exceeding $x$.
 Let $q$ be a
prime of the form $4k+1$ and so $q=a^2+b^2$ with $a,b\in\Bbb Z$.
When $p\nmid ab(a^2-b^2)q$, we show that for $r=0,1,2,3$,
$p^{\frac{q-1}4}\equiv (\frac ab)^r\pmod q$ if and only if
$$\sum_{k=0}^{[p/4]}\binom{4k}{2k}\Big(\frac{a^2}{16q}\Big)^k\equiv
(-1)^{\frac{p^2-1}8a+\frac{p-1}2\cdot \frac{q-1}4}\Big(\frac pq\Big)
\Big(\frac ab\Big)^r\pmod p,$$
 where $(\frac pq)$ is the Legendre symbol.
 We also establish congruences for
$\sum_{k=0}^{[p/4]}\binom{4k}{2k}\frac 1{m^k}\pmod p$ in the cases
$m=17,18,20,32,52,80,272$.

\par\q
 }
 \endabstract
 \footnotetext[1] {The author is
supported by the Natural Sciences Foundation of China (grant no.
11371163).}

\section*{1. Introduction}
\par\q\par Let $\Bbb Z$ be the set of integers, and for a prime
$p$ let $\Bbb Z_p$ denote the set of  those rational numbers whose
denominator is not divisible by $p$. Let $\sls mp$ be the Legendre
symbol.

\par Suppose that $p$ is an odd prime and $a\in\Bbb Z_p$.  In
[7] the author investigated congruences for
$\sum_{k=0}^{[p/4]}\b{4k}{2k}a^k$ modulo $p$, where $[x]$ is the
greatest integer not exceeding $x$.  For $k\in\{0,1,\ldots,p-1\}$ it
is easily seen that $p\nmid \b{4k}{2k}$ if and only if $0\le k<\f
p4$ or $\f p2<k<\f{3p}4$. In this paper we reveal the connection
between quartic residues and the sum
$\sum_{k=0}^{[p/4]}\b{4k}{2k}a^k$. We also investigate congruences
for $\sum_{k=(p+1)/2}^{[3p/4]}\b{4k}{2k}a^k$ modulo $p$.
\par Let  $i=\sqrt{-1}$.  For an odd prime $p$ let $\ls{a+bi}p_4$
be the quartic Jacobi symbol defined in [1,2,3,4,6]. Following [4]
we define
$$Q_r(p)=\Big\{c\in\Bbb Z_p: \Ls{c+i}p_4=i^r\Big\}\qtq{for}
r=0,1,2,3.$$
 According
to [4,6], $Q_r(p)\ (r=0,1,2,3)$ play a central role in the theory of
quartic residues and nonresidues. In this paper, for an odd prime
$p$ and $c\in\Bbb Z_p$ with $c(c^2+1)\not\e 0\mod p$ we give simple
criteria for $c\in Q_r(p)$ by proving that
$$\sum_{k=0}^{[p/4]}\b{4k}{2k}\f 1{(16(c^2+1))^k}
\e\cases 1\mod p&\t{if $c\in Q_0(p)$,}
\\-\f 1c\mod p&\t{if $c\in Q_1(p)$,}
\\-1\mod p&\t{if $c\in Q_2(p)$,}
\\\f 1c\mod p&\t{if $c\in Q_3(p)$}\endcases\tag 1.1$$
and
$$ 2c\sum_{k=(p+1)/2}^{[3p/4]}\b{4k}{2k}\f 1{(16(c^2+1))^k}
\e \cases 0\mod p&\t{if $\sls{c^2+1}p=1$,}
\\1\mod p&\t{if $c\in Q_1(p)$,}
\\-1\mod p&\t{if $c\in Q_3(p)$.}
\endcases\tag 1.2$$
\par Let $q$ be a prime of the form $4m+1$ and so $q=a^2+b^2$ with
$a,b\in\Bbb Z$. Let $p$ be an odd prime with $p\nmid ab(a^2-b^2)q$.
In this paper, using (1.1) and the theory of quartic residues we
show that for $r=0,1,2,3$,
$$\aligned &p^{\f{q-1}4}\e \Ls ab^r\mod q\\&\Leftrightarrow
\sum_{k=0}^{[p/4]}\b{4k}{2k}\Ls{a^2}{16q}^k\e
(-1)^{\f{p^2-1}8a+\f{p-1}2\cdot \f{q-1}4}\Ls pq\Ls ab^r\mod
p.\endaligned\tag 1.3$$
  As
consequences we establish congruences for
$\sum_{k=0}^{[p/4]}\b{4k}{2k}\f 1{m^k}\mod p$ in the cases
$m=17,18,20,32,52,80,272$.
\par In addition to the above notation, throughout this paper we use
$(m,n)$ to denote the greatest common divisor of integers $m$ and
$n$. If $a,b,c\in\Bbb Z$ and $p=ax^2+bxy+cy^2$ for some integers $x$
and $y$, we briefly write that $p=ax^2+bxy+cy^2$. We also use
$[a,b,c]$ to denote the equivalence class containing the form
$ax^2+bxy+cy^2$, and use $H(d)$ to denote the form class group
consisting of equivalence classes of discriminant $d$.
\section*{2. Congruences for $\sum_{k=0}^{[p/4]}\b{4k}{2k}
a^k\mod p$}

\par For any numbers $P$ and $Q$, let $\{U_n(P,Q)\}$ and
$\{V_n(P,Q)\}$ be the Lucas sequences given by
$$\align &U_0(P,Q)=0,\ U_1(P,Q)=1,\ U_{n+1}(P,Q)=PU_n(P,Q)-QU_{n-1}(P,Q)
\ (n\ge 1),\\&V_0(P,Q)=2,\ V_1(P,Q)=P,\
V_{n+1}(P,Q)=PV_n(P,Q)-QV_{n-1}(P,Q) \ (n\ge 1).\endalign$$ It is
well known (see [8]) that
$$\aligned &U_n(P,Q)\\&=\cases \f
1{\sqrt{P^2-4Q}}\Big\{\Big(\f{P+\sqrt{P^2-4Q}}2\Big)^n-
\Big(\f{P-\sqrt{P^2-4Q}}2\Big)^n\Big\} &\t{if $P^2-4Q\not=0$,}
\\n\ls P2^{n-1}&\t{if $P^2-4Q=0$.}
\endcases\endaligned\tag 2.1$$
and
$$V_n(P,Q)=\Big(\f{P+\sqrt{P^2-4Q}}2\Big)^n+
\Big(\f{P-\sqrt{P^2-4Q}}2\Big)^n.\tag 2.2$$ \pro{Lemma 2.1 ([7,
Theorem 2.1])} Let $p$ be an odd prime, $P,Q\in\Bbb Z_p$ and
$PQ\not\e 0\mod p$. Then
$$\sum_{k=0}^{[p/4]}\b{4k}{2k}\Ls Q{4P^2}^k\e \Ls{P}p
U_{\f{p+1}2}(P,Q)\mod p$$ and
$$\sum_{k=0}^{[p/4]}\b{4k}{2k}\Big(\f{P^2}{64Q}\Big)^k\e
(-Q)^{-[\f p4]}U_{\f{p+\sls{-1}p}2}(P,Q)\mod p.$$
\endpro

\pro{Lemma 2.2 ([7, Theorem 2.2])} Let $p$ be an odd prime and
$x\in\Bbb Z_p$ with $x\not\e 0,1\mod p$. Then
$$\sum_{k=0}^{[p/4]}\b{4k}{2k}
\f 1{(16x)^k}\e\cases \f 1
{x^{\f{p-1}4}}\sum_{k=0}^{[p/4]}\b{4k}{2k}\Ls x{16}^k \mod p&\t{if
$4\mid p-1$,}
\\ \Big(1-\f 1x\Big)^{\f{p-3}4}\sum_{k=0}^{[p/4]}\b{4k}{2k}
\f 1{(16(1-x))^k}\mod p&\t{if $4\mid p-3$}.\endcases$$
\endpro
\pro{Lemma 2.3 ([4, Lemma 2.1])} Let $p$ be an odd prime,
$m,n\in\Bbb Z_p$ and $m^2+n^2\not\e 0\mod p$. Then
$\ls{m+ni}p_4^2=\ls{m^2+n^2}p$.
\endpro
 \pro{Theorem 2.1} Suppose that $p$ is an odd prime, $c\in\Bbb Z_p$
 and $c(c^2+1)\not\e 0\mod p$. Then
$$\Ls 2p\sum_{k=0}^{[p/4]}\b{4k}{2k}\Ls{c^2}{16(c^2+1)}^k
\e\cases 1\mod p&\t{if $c\in Q_0(p)$,}
\\c\mod p&\t{if $c\in Q_1(p)$,}
\\-1\mod p&\t{if $c\in Q_2(p)$,}
\\-c\mod p&\t{if $c\in Q_3(p)$}\endcases\tag 2.3$$
and
$$\sum_{k=0}^{[p/4]}\b{4k}{2k}\f 1{(16(c^2+1))^k}
\e\cases 1\mod p&\t{if $c\in Q_0(p)$,}
\\-\f 1c\mod p&\t{if $c\in Q_1(p)$,}
\\-1\mod p&\t{if $c\in Q_2(p)$,}
\\\f 1c\mod p&\t{if $c\in Q_3(p)$.}\endcases\tag 2.4$$
\endpro
Proof. Clearly
$$\Ls{c+i}p_4=\Ls ip_4\Ls{1-ci}p_4=\Ls 2p\Ls{-\f 1c+i}p_4.$$
Thus, if $\ls 2p=1$, then $-\f 1c\in Q_r(p)$ if and only if $c\in
Q_r(p)$; if $\ls 2p=-1$, then $-\f 1c\in Q_r(p)$ if and only if
$c\in Q_{r'}(p)$, where $r'\in\{0,1,2,3\}$ is given by $r'\e r+2\mod
4$. Thus, replacing $c$ with $-\f 1c$ in (2.3) we get (2.4). Hence
(2.3) is equivalent to (2.4).
\par By [7, proof of Theorem
2.1], for $P,Q\in\Bbb Z_p$ with $PQ\not\e 0\mod p$,
$$V_{\f{p-1}2}(P,(P^2-4Q)/4)
\e 2\Ls {2P}p\sum_{k=0}^{[p/4]}\b{4k}{2k}\Ls Q{4P^2}^k\mod p.\tag
2.5$$ Taking $P=2c$ and $Q=c^2+1$ we see that
$$V_{\f{p-1}2}(2c,-1)\e 2\Ls
cp\sum_{k=0}^{[p/4]}\b{4k}{2k}\Ls{c^2+1}{16c^2}^k\mod p.$$ From
Lemma 2.3 we have $\ls {n+i}p_4^2=\ls{n^2+1}p$ for $n^2+1\not\e
0\mod p$. We first assume $p\e 1\mod 4$. Taking $a=-1$ and $b=2c$ in
[5, Corollary 3.1(i)] we get
$$V_{\f{p-1}2}(2c,-1)\e \cases 2(4c^2+4)^{\f{p-1}4}\mod p&\t{if $c\in
Q_0(p)$,}\\2c(4c^2+4)^{\f{p-1}4}\mod p&\t{if $c\in Q_1(p)$,}
\\-2(4c^2+4)^{\f{p-1}4}\mod p&\t{if $c\in
Q_2(p)$,}\\-2c(4c^2+4)^{\f{p-1}4}\mod p&\t{if $c\in Q_3(p)$.}
\endcases$$
Combining  the above with the fact $4^{\f{p-1}4}=2^{\f{p-1}2}\e \sls
2p\mod p$ we obtain
$$\Ls {2c}p\sum_{k=0}^{[p/4]}\b{4k}{2k}\Ls{c^2+1}{16c^2}^k
\e\cases (c^2+1)^{\f{p-1}4}\mod p&\t{if $c\in Q_0(p)$,}
\\c(c^2+1)^{\f{p-1}4}\mod p&\t{if $c\in Q_1(p)$,}
\\-(c^2+1)^{\f{p-1}4}\mod p&\t{if $c\in Q_2(p)$,}
\\-c(c^2+1)^{\f{p-1}4}\mod p&\t{if $c\in Q_3(p)$.}\endcases
\tag 2.6$$ Putting $x=\f{c^2+1}{c^2}$ in Lemma 2.2 we see that
$$\sum_{k=0}^{[p/4]}\b{4k}{2k}\Ls{c^2}{16(c^2+1)}^k\e
(c^2+1)^{-\f{p-1}4}\Ls cp
\sum_{k=0}^{[p/4]}\b{4k}{2k}\Ls{c^2+1}{16c^2}^k\mod p.$$ This
together with (2.6) yields (2.3). Hence (2.4) is also true.
\par Now we assume $p\e 3\mod 4$.
Taking $P=2c$ and $Q=-1$ in Lemma 2.1 we see that
$$\sum_{k=0}^{[p/4]}\b{4k}{2k}\f 1{(-16c^2)^k}
\e \Ls{2c}pU_{\f{p+1}2}(2c,-1)\mod p.$$ By [5, Theorem 3.1(ii)] and
Lemma 2.3,
$$U_{\f{p+1}2}(2c,-1)
\e\cases 2c(4c^2+4)^{\f{p-3}4}\mod p&\t{if $c\in Q_0(p)$,}
\\-2(4c^2+4)^{\f{p-3}4}\mod p&\t{if $c\in Q_1(p)$,}
\\-2c(4c^2+4)^{\f{p-3}4}\mod p&\t{if $c\in Q_2(p)$,}
\\2(4c^2+4)^{\f{p-3}4}\mod p&\t{if $c\in Q_3(p)$.}\endcases$$
Thus,
$$\Ls cp\sum_{k=0}^{[p/4]}\b{4k}{2k}\f 1{(-16c^2)^k}
\e\cases c(c^2+1)^{\f{p-3}4}\mod p&\t{if $c\in Q_0(p)$,}
\\-(c^2+1)^{\f{p-3}4}\mod p&\t{if $c\in Q_1(p)$,}
\\-c(c^2+1)^{\f{p-3}4}\mod p&\t{if $c\in Q_2(p)$,}
\\(c^2+1)^{\f{p-3}4}\mod p&\t{if $c\in Q_3(p)$.}\endcases$$
By Lemma 2.2,
$$\sum_{k=0}^{[p/4]}\b{4k}{2k}\f 1{(16(c^2+1)^k}
\e \f 1c\Ls cp(c^2+1)^{-\f{p-3}4}\sum_{k=0}^{[p/4]}\b{4k}{2k}\f
1{(-16c^2)^k}\mod p.$$ Now combining the above we obtain (2.4).
Since (2.3) is equivalent to (2.4), the proof is complete.
\par\q
\par{\bf Remark 2.1} Let $p$ be an odd prime, $a\in\Bbb Z_p$,
$a\not\e -1\mod p$ and $\sls ap=-1$. Taking $P=2$ and $Q=a+1$ in [7,
Theorem 2.3(i)] we see that
$$\sum_{k=0}^{[p/4]}\b{4k}{2k}\f 1{(16(a+1))^k}\e 0\mod p.$$
Replacing $a$ with $\f 1a$ we get
$$\sum_{k=0}^{[p/4]}\b{4k}{2k}\Ls a{16(a+1)}^k\e 0\mod p.$$
 \pro{Corollary 2.1} Suppose that $p$ is an odd prime, $c\in\Bbb
Z_p$
 and $c(c^2+1)\not\e 0\mod p$. Then
$$\Big(\sum_{k=0}^{[p/4]}\b{4k}{2k}\Ls{c^2}{16(c^2+1)}^k\Big)
\Big(\sum_{k=0}^{[p/4]}\b{4k}{2k}\f 1{(16(c^2+1))^k}\Big) \e
\Ls{2(c^2+1)}p\mod p.$$
\endpro
Proof. By Lemma 2.3, $c\in Q_0(p)\cup Q_2(p)$ if and only if
$\ls{c^2+1}p=1$. Thus the result follows from Theorem 2.1.
\pro{Corollary 2.2} Suppose that $p$ and $q$ are distinct primes,
$m,n\in\Bbb Z$, $(mn(m^2+n^2)(m^2-n^2),pq)=1$ and $p\e \pm
q\mod{8(m^2+n^2)}$. Assume $a\in\{1,-1,\f nm,-\f nm\}$. Then
$$\align &\sum_{k=0}^{[p/4]}\b{4k}{2k}\Ls{n^2}{16(m^2+n^2)}^k\e
a\mod p
\\&\iff \sum_{k=0}^{[q/4]}\b{4k}{2k}\Ls{n^2}{16(m^2+n^2)}^k\e
a\mod q.\endalign$$
\endpro
Proof. By [4, Theorem 2.1], $\f mn\in Q_r(p)\Leftrightarrow \f mn\in
Q_r(q)$. Now taking $c=\f mn$ in (2.4) we obtain the result.

\pro{Corollary 2.3} Let $p$ be an odd prime. Then
$$\sum_{k=0}^{[p/4]}\b{4k}{2k}\f 1{32^k}
\e\cases 1\mod p&\t{if $p\e \pm 1,\pm 3\mod {16}$,} \\-1\mod p&\t{if
$p\e \pm 5,\pm 7\mod {16}$.}\endcases$$
\endpro
Proof. It is well known (see [1-3]) that
$$\Ls{1+i}p_4=i^{\f{(-1)^{\f{p-1}2}p-1}4}=\cases
1&\t{if $p\e \pm 1\mod{16}$,}
\\i&\t{if $p\e \pm 5\mod{16}$,}
\\-1&\t{if $p\e \pm 7\mod{16}$,}
\\-i&\t{if $p\e \pm 3\mod{16}$.}
\endcases\tag 2.7$$
Now taking $c=1$ in (2.4) we deduce the result.

\pro{Corollary 2.4} Suppose that $p$ is an odd prime, $c\in\Bbb
Z_p$, $c\not\e 0,\pm 1\mod p$ and $c^2+1\not\e 0\mod p$. Then the
congruence $x^4-2(c^2+1)x^2+c^2(c^2+1)\e 0\mod p$ is solvable if and
only if
$$\sum_{k=1}^{[p/4]}\b{4k}{2k}\f 1{(16(c^2+1))^k}
\e 0\mod p.$$\endpro Proof. By [6, Theorem 4.2], $c\in C_0(p)$ if
and only if $x^4-2(c^2+1)x^2+c^2(c^2+1)\e 0\mod p$ is solvable. By
Theorem 2.1, $c\in C_0(p)$ if and only if
$\sum_{k=0}^{[p/4]}\b{4k}{2k}\f 1{(16(c^2+1))^k} \e 1\mod p.$ Thus
the result follows.

\pro{Corollary 2.5} Let $p$ be an odd prime, $d\in\Bbb Z$, $d\not\e
0,\pm 1\mod p$, $(\f{-d}p)=1$, $n\in\Bbb Z_p$, $n^2\e -d\mod p$ and
$1-d=(-1)^r2^sW\;(W\e 1\mod 4)$. Let $d_0$ be the product of all
distinct odd prime divisors of $1-d$, and let $m=8d_0$ or
$\f{4d_0}{(2,r+s/2)}$ according as $2\nmid s$ or $2\mid s$.
 Then  $p$ is represented by some primitive form $ax^2+2bxy+cy^2$
of discriminant $-4m^2d$ with the condition that $(a,2(1-d))=1$ and
$\qs{b-mi}a=1$ if and only if
$$\sum_{k=1}^{[p/4]}\b{4k}{2k}\f 1{(16(1-d))^k}
\e 0\mod p.$$\endpro Proof. This is immediate from [6, Theorem 4.2]
and (2.4).
\par It is known that
$$\align &H(-128)=\{[1,0,32],[4,4,9],[3,2,11],[3,-2,11]\},
\\&H(-768)=\{[1,0,192],[12,12,19],[3,0,64],[4,4,49],
\\&\qq\qq\qq [7,4,28],[7,-4,28],
[13,8,16],[13,-8,16]\},\\&
H(-80)=\{[1,0,20],[4,0,5],[3,2,7],[3,-2,7]\}.\endalign$$
 Thus taking $d=2,3,5$ in Corollary 2.5 and doing some calculations for
 certain quartic Jacobi symbols we see
that for any prime $p>3$,
$$\align &p=x^2+32y^2\Leftrightarrow \sum_{k=1}^{[p/4]}\b{4k}{2k}\f 1{(-16)^k}
\e 0\mod p,\tag 2.8\\&p=x^2+192y^2\ \t{or}\ 12x^2+12xy+19y^2
\Leftrightarrow \sum_{k=1}^{[p/4]}\b{4k}{2k}\f 1{(-32)^k} \e 0\mod
p,\tag 2.9\\&p=x^2+20y^2\Leftrightarrow
\sum_{k=1}^{[p/4]}\b{4k}{2k}\f 1{(-64)^k} \e 0\mod p.\tag
2.10\endalign$$

 \pro{Lemma 2.4 ([4, Theorem 2.2])} Let $q$ be a prime of
 the form $4m+1$ and
$q=a^2+b^2$ with $a,b\in\Bbb Z$, $2\mid b$ and $a+b\e 1\mod 4$.
Suppose that $p$ is an odd prime and $p\nmid bq$. For $r=0,1,2,3$ we
have
$$((-1)^{\f{p-1}2}p)^{\f{q-1}4}\e \big(-\f ab\big)^r\mod q
\iff \f ab\in Q_r(p).$$
\endpro

\pro{Theorem 2.2} Let $q$ be a prime of the form $4m+1$ and
$q=a^2+b^2$ with $a,b\in\Bbb Z$ and $a\e 1\mod 2$. Suppose that $p$
is an odd prime and $p\nmid bq$. Then
$$(-1)^{\f{p^2-1}8+\f{p-1}2\cdot \f{q-1}4}
\sum_{k=0}^{[p/4]}\b{4k}{2k}\Ls{a^2}{16q}^k\e \cases \pm 1\mod
p&\t{if $p^{\f{q-1}4}\e\pm 1\mod q$,}
\\\mp \f ab\mod p&\t{if $p^{\f{q-1}4}\e\pm \f ab\mod q$.}
\endcases$$
\endpro
Proof. Clearly the result does not depend on the sign of $a$. We may
assume $a+b\e 1\mod 4$. Taking $c=\f ab$ in (2.3) we see that
$$\Ls 2p\sum_{k=0}^{[p/4]}\b{4k}{2k}\Ls{a^2}{16q}^k
\e\cases 1\mod p&\t{if $\f ab\in Q_0(p)$,}
\\\f ab\mod p&\t{if $\f ab\in Q_1(p)$,}
\\-1\mod p&\t{if $\f ab\in Q_2(p)$,}
\\-\f ab\mod p&\t{if $\f ab\in Q_3(p)$.}\endcases$$
This is also true when $p\mid a$ since $\sls{0+i}p_4=\sls 2p$. Now
applying Lemma 2.4 we deduce the result.
 \pro{Corollary 2.6} Let
$p\not=5$ be an odd prime. Then
$$(-1)^{[\f p4]}\sum_{k=0}^{[p/4]}\b{4k}{2k}\f 1{80^k}
\e\cases \pm 1\mod p&\t{if $p\e \pm 1\mod 5$,}
\\\pm \f 12\mod p&\t{if $p\e \pm 2\mod 5$.}
\endcases$$
\endpro
Proof. Taking $q=5$, $a=1$ and $b=2$ in Theorem 2.2 we deduce the
result.

\pro{Corollary 2.7} Let $p\not=13$ be an odd prime. Then
$$(-1)^{[\f p4]}\sum_{k=0}^{[p/4]}\b{4k}{2k}\Ls 9{208}^k
\e\cases \pm 1\mod p&\t{if $p\e \pm 1,\pm 3,\pm 9\mod {13}$,}
\\\mp \f 32\mod p&\t{if $p\e \pm 2,\pm 5,\pm 6\mod {13}$.}
\endcases$$
\endpro
Proof. Taking $q=13$, $a=-3$ and $b=2$ in Theorem 2.2 we deduce the
result.

\pro{Corollary 2.8} Let $p\not=17$ be an odd prime. Then
$$\Ls 2p\sum_{k=0}^{[p/4]}\b{4k}{2k}\f 1{272^k}
\e\cases  1\mod p&\t{if $p\e \pm 1,\pm 4\mod {17}$,}
\\-1\mod p&\t{if $p\e \pm 2,\pm 8\mod {17}$,}
\\ \f 14\mod p&\t{if $p\e \pm 6,\pm 7\mod {17}$,}
\\-\f 14\mod p&\t{if $p\e \pm 3,\pm 5\mod {17}$.}
\endcases$$
\endpro
Proof. Taking $q=17$, $a=1$ and $b=4$ in Theorem 2.2 we deduce the
result.

\pro{Theorem 2.3} Let $q$ be a prime of the form $4m+1$ and
$q=a^2+b^2$ with $a,b\in\Bbb Z$ and $a\e 1\mod 2$. Suppose that $p$
is an odd prime and $p\nmid aq$. Then
$$(-1)^{\f{p-1}2\cdot \f{q-1}4}
\sum_{k=0}^{[p/4]}\b{4k}{2k}\Ls{b^2}{16q}^k \e\cases \pm 1\mod
p&\t{if $p^{\f{q-1}4}\e\pm 1\mod q$,}
\\\mp \f ba\mod p&\t{if $p^{\f{q-1}4}\e\pm \f ba\mod q$.}
\endcases$$
\endpro
Proof. Clearly the result does not depend on the sign of $a$. We may
assume $a+b\e 1\mod 4$. When $p\mid b$, by [4, Theorem 2.2] we have
$((-1)^{\f{p-1}2}p)^{\f{q-1}4}\e 1\mod q$. Thus the result is true
in this case. Now we assume $p\nmid b$. Taking $c=\f ab$ in (2.4)
and then applying Lemma 2.4 we obtain the result.
\par\q
\par Combining Theorem 2.2 with Theorem 2.3 we obtain (1.3).
 \pro{Corollary 2.9} Let $p$ be an odd prime. Then
$$\align&(-1)^{\f{p-1}2}\sum_{k=0}^{[p/4]}\b{4k}{2k}\f 1{20^k}
\e\cases \pm 1\mod p&\t{if $p\e \pm 1\mod 5$,}
\\\pm 2\mod p&\t{if $p\e \pm 3\mod 5$,}\endcases
\\&(-1)^{\f{p-1}2}\sum_{k=0}^{[p/4]}\b{4k}{2k}\f 1{52^k}
\e\cases \pm 1\mod p&\t{if $p\e \pm 1,\pm 3,\pm 9\mod {13}$,}
\\\pm \f 23\mod p&\t{if $p\e \pm 2,\pm 5,\pm 6\mod {13}$.}
\endcases\endalign$$
and
$$\sum_{k=0}^{[p/4]}\b{4k}{2k}\f 1{17^k}
\e\cases  1\mod p&\t{if $p\e \pm 1,\pm 4\mod {17}$,}
\\-1\mod p&\t{if $p\e \pm 2,\pm 8\mod {17}$,}
\\ 4\mod p&\t{if $p\e \pm 3,\pm 5\mod {17}$,}
\\-4\mod p&\t{if $p\e \pm 6,\pm 7\mod {17}$.}
\endcases$$
\endpro
\endpro
Proof. Taking $q=5,13,17$ in Theorem 2.3 we deduce the result.

\pro{Theorem 2.4} Let $p$ be an odd prime, $a,b\in\Bbb Z_p$ and
$ab(a^2-b^2)\not\e 0\mod p$. Then
$$\sum_{k=0}^{[p/4]}\b{4k}{2k}\Ls {a^2-b^2}{16a^2}^k
\e\f 1{2b}\Ls{2a}p\Big\{(a+b)\Ls {a+b}p-(a-b)\Ls{a-b}p\Big\}\mod p$$
and
$$\align&\sum_{k=0}^{[p/4]}\b{4k}{2k}\Ls{a^2}{16(a^2-b^2)}^k
\\&\e \cases\f 1{2b}\big\{(a+b)\ls
{a-b}p-(a-b)\ls{a+b}p\big\}(b^2-a^2)^{\f{p-1}4}\mod p &\t{if $4\mid
p-1$,}
\\\f 1{2b}\big\{\ls
{a+b}p-\ls{a-b}p\big\}(b^2-a^2)^{\f{p+1}4}\mod p &\t{if $4\mid
p-3$.}\endcases
\endalign$$
\endpro
Proof.  By (2.1) we have
$$\align U_{\f{p+1}2}\Big(a,\f{a^2-b^2}4\Big)&=\f
1b\Big\{\Ls{a+b}2^{\f{p+1}2}-\Ls{a-b}2^{\f{p+1}2}\Big\}
\\&\e \f 1{2b}\Ls 2p
\Big\{(a+b)\Ls{a+b}p-(a-b)\Ls{a-b}p\Big\}\mod p.\endalign$$ Now
taking $P=a$ and $Q=\f{a^2-b^2}4$ in Lemma 2.1 and then applying the
above we obtain the first part.
\par Taking $P=2a$ and $Q=a^2-b^2$ in Lemma 2.1 and then applying (2.1)
we see that
$$\align &\sum_{k=0}^{[p/4]}\b{4k}{2k}\Ls{a^2}{16(a^2-b^2)}^k
\\&\e (b^2-a^2)^{-[\f p4]}U_{\f{p+\sls{-1}p}2}(2a,a^2-b^2)
\\&=(b^2-a^2)^{\f{p-\sls{-1}p}4-\f{p-1}2}\f
1{2b}\Big\{(a+b)^{\f{p+\sls{-1}p}2}-(a-b)^{\f{p+\sls{-1}p}2}\Big\}
\\&\e\cases
\ls{b^2-a^2}p(b^2-a^2)^{\f{p-1}4} \f
1{2b}\big\{(a+b)\ls{a+b}p-(a-b)\ls{a-b}p\big\} \mod p&\t{if $4\mid
p-1$,}
\\\ls{b^2-a^2}p(b^2-a^2)^{\f{p+1}4} \f
1{2b}\big\{\ls{a+b}p-\ls{a-b}p\big\} \mod p&\t{if $4\mid p-3$.}
\endcases\endalign$$
This yields the remaining part.
 \pro{Corollary 2.10} Let $p$ be an odd
prime, $m\in\Bbb Z_p$ and $m\not\e 0,\pm 1\mod p$. Then
$$\sum_{k=0}^{[p/4]}\b{4k}{2k}\Big(-\f m{4(m-1)^2}\Big)^k
\e\f 1{m+1}\Ls{m-1}p\Big\{m\Ls mp+\Ls{-1}p\Big\}\mod p$$ and
$$\sum_{k=0}^{[p/4]}\b{4k}{2k}\Big(-\f {(m-1)^2}{64m}\Big)^k
\e \cases\f 1{m+1}\big(m+\ls mp\big)m^{\f{p-1}4}\mod p&\t{if $4\mid
p-1$,}\\\f 1{m+1}\big(\ls mp+1\big)m^{\f{p+1}4}\mod p&\t{if $4\mid
p-3$.}
\endcases$$
\endpro
Proof. Taking $a=\f{m-1}2$ and $b=\f{m+1}2$ in Theorem 2.4 we deduce
the result.
 \pro{Corollary 2.11} Let $p>3$ be a
prime. Then
$$\sum_{k=0}^{[p/4]}\b{4k}{2k}\f 1{18^k}
\e 2\Ls 6p-\Ls 3p\mod p.$$
\endpro
Proof. Taking $a=3$ and $b=1$ in Theorem 2.4 we obtain the result.
\section*{3. Congruences for $\sum_{k=(p+1)/2}^{[3p/4]}\b{4k}{2k}
a^k\mod p$}
\par
 \pro{Lemma 3.1} Let $p$ be an odd prime
and $k\in\{1,2,\ldots,\f {p-1}2\}$. Then
$$\b{4(\f{p-1}2+k)}{2(\f{p-1}2+k)}\e 2\b{4k-2}{2k-1}\mod p.$$
\endpro
Proof. Clearly
$$\align \b{4(\f{p-1}2+k)}{2(\f{p-1}2+k)}&=
\f{(p+2k)(p+2k+1)\cdots(2p-1)\cdot 2p(2p+1)\cdots (2p+4k-2)}
{(p-1)!\cdot p(p+1)\cdots(p+2k-1)}
\\&\e \f{2k(2k+1)\cdots (p-1)\cdot 2\cdot 1\cdot 2\cdots (4k-2)}
{(p-1)!\cdot 1\cdot 2\cdots (2k-1)}
\\&=\f{2\cdot (4k-2)!}{(2k-1)!^2}=2\b{4k-2}{2k-1}\mod p.
\endalign$$
\pro{Lemma 3.2} Let $p$ be an odd prime, $P,Q\in\Bbb Z_p$ and
$PQ\not\e 0\mod p$. Then
$$U_{\f{p-1}2}(P,Q)\e \f{2P}Q\Ls{P}p\sum_{k=1}^{[\f{p+1}4]}\b{4k-2}{2k-1}\Ls
Q{4P^2}^k\mod p.$$\endpro

Proof. By (2.1) and the fact $\b{\f{p-1}2}k\e \b{-\f 12}k =\f
1{(-4)^k}\b{2k}k\mod p$ for $k\in\{0,1,\ldots,\f{p-1}2\}$,
$$\align U_{\f{p-1}2}(P,(P^2-4Q)/4)
&=\f 1{2\sqrt Q}\Big\{ \Ls{P+2\sqrt Q}2^{\f{p-1}2}-\Ls{P-2\sqrt
Q}2^{\f{p-1}2}\Big\}
\\&=\f 1{2^{\f{p-1}2}\cdot 2\sqrt Q}
\sum_{r=0}^{\f{p-1}2}\b{\f{p-1}2}rP^{\f{p-1}2-r}\big((2\sqrt Q)^r-
(-2\sqrt Q)^r\big)
\\&=\f 1{2^{\f{p-1}2}\cdot 2\sqrt Q}
\sum_{k=1}^{[\f{p+1}4]}\b{\f{p-1}2}{2k-1}P^{\f{p-1}2-(2k-1)}\cdot
2(2\sqrt Q)^{2k-1}
\\&\e \Ls 2p\sum_{k=1}^{[\f{p+1}4]}\b{4k-2}{2k-1}\f 1{(-4P)^{2k-1}}
\Ls Pp\cdot 2(4Q)^{k-1}
\\&=-\f{2P}Q\Ls{2P}p\sum_{k=1}^{[\f{p+1}4]}\b{4k-2}{2k-1}\Ls
Q{4P^2}^k\mod p.\endalign$$ By [5, Lemma 3.1(i)], $U_{\f{p-1}2}(P,Q)
\e -\sls 2pU_{\f{p-1}2}(P,\f{P^2-4Q}4)\mod p$. Thus the result
follows.

\pro{Theorem 3.1} Let $p$ be an odd prime, $c\in\Bbb Z_p$ and
$c(c^2+1)\not\e 0\mod p$. Then
$$\align 2c\sum_{k=(p+1)/2}^{[3p/4]}\b{4k}{2k}\f 1{(16(c^2+1))^k}
&\e -4c\sum_{k=1}^{[\f{p+1}4]}\b{4k-2}{2k-1} \f
1{(16(c^2+1))^k}\\&\e \cases 0\mod p&\t{if $\sls{c^2+1}p=1$,}
\\1\mod p&\t{if $c\in Q_1(p)$,}
\\-1\mod p&\t{if $c\in Q_3(p)$.}
\endcases\endalign$$
\endpro
Proof. By Lemmas 3.1 and 3.2,
$$\align
\sum_{k=(p+1)/2}^{[3p/4]}\b{4k}{2k}\f 1{(16(c^2+1))^k}
&=\sum_{k=1}^{[\f{p+1}4]}\b{4(\f{p-1}2+k)}{2(\f{p-1}2+k)} \f
1{(16(c^2+1))^{\f{p-1}2+k}}
\\&\e 2\Ls{c^2+1}p\sum_{k=1}^{[\f{p+1}4]}\b{4k-2}{2k-1}
\f 1{(16(c^2+1))^k}
\\&\e U_{\f{p-1}2}(4(c^2+1),4(c^2+1))\mod p.\endalign$$
From (2.1) we see that
$$\align &U_{\f{p-1}2}(4(c^2+1),4(c^2+1))\\&
=\f 1{4c\sqrt{c^2+1}}(2\sqrt{c^2+1})^{\f{p-1}2}
\big\{(\sqrt{c^2+1}+c)^{\f{p-1}2}-(\sqrt{c^2+1}-c)^{\f{p-1}2}\big\}
\\&=\cases\f 1{2c}(4c^2+4)^{\f{p-1}4}U_{\f{p-1}2}(2c,-1)&\t{if
$4\mid p-1$,}
\\\f 1{2c}(4c^2+4)^{\f{p-3}4}V_{\f{p-1}2}(2c,-1)
&\t{if $4\mid p-3$.}\endcases\endalign$$ If $p\e 1\mod 4$, by [5,
Theorem 3.1(i)] we have
$$U_{\f{p-1}2}(2c,-1)
\e\cases 0\mod p&\t{if $\sls{c^2+1}p=1$,}
\\(4c^2+4)^{\f{p-1}4}\qs{c+i}pi\mod p&\t{if $\sls{c^2+1}p=-1$.}
\endcases$$
Thus,
$$\align 2c\sum_{k=(p+1)/2}^{[3p/4]}\b{4k}{2k}\f 1{(16(c^2+1))^k}
&\e (4c^2+4)^{\f{p-1}4}U_{\f{p-1}2}(2c,-1)
\\&\e\cases 0\mod p &\t{if
$\sls{c^2+1}p=1$,}
\\-\qs{c+i}pi\mod p&\t{if $\sls{c^2+1}p=-1$.}
\endcases\endalign$$
If $p\e 3\mod 4$, from [5, Corollary 3.1(ii)] we have
$$V_{\f{p-1}2}(2c,-1)
\e\cases 0\mod p&\t{if $\sls{c^2+1}p=1$,}
\\(4c^2+4)^{\f{p+1}4}\qs{c+i}pi\mod p&\t{if $\sls{c^2+1}p=-1$.}
\endcases$$
Thus,
$$\align 2c\sum_{k=(p+1)/2}^{[3p/4]}\b{4k}{2k}\f 1{(16(c^2+1))^k}
&\e (4c^2+4)^{\f{p-3}4}V_{\f{p-1}2}(2c,-1)
\\&\e\cases 0\mod p &\t{if
$\sls{c^2+1}p=1$,}
\\-\qs{c+i}pi\mod p&\t{if $\sls{c^2+1}p=-1$.}
\endcases\endalign$$
Note that $\qs{c+i}p^2=\sls{c^2+1}p$ by Lemma 2.3.
Combining all the
above we deduce the result.

\pro{Corollary 3.1} Let $p$ be an odd prime, $c\in\Bbb Z_p$ and
$c(c^2+1)\not\e 0\mod p$. Then
$$\align -\f 2c\Ls 2p\sum_{k=(p+1)/2}^{[3p/4]}\b{4k}{2k}\Ls {c^2}
{16(c^2+1)}^k &\e \f 4c\Ls 2p\sum_{k=1}^{[\f{p+1}4]}\b{4k-2}{2k-1}
\Ls {c^2}{16(c^2+1)}^k
\\&\e \cases 0\mod p&\t{if $\sls{c^2+1}p=1$,}
\\1\mod p&\t{if $c\in Q_1(p)$,}
\\-1\mod p&\t{if $c\in Q_3(p)$.}
\endcases\endalign$$
\endpro
Proof. Clearly $\ls{c+i}p_4=\ls ip_4\ls{1-ci}p_4=\ls 2p\ls{-\f
1c+i}p_4.$ If $\sls 2p=1$, then $-\f 1c\in Q_r(p)$ if and only if
$c\in Q_r(p)$. If $\sls 2p=-1$, then $-\f 1c\in Q_1(p)$ if and only
if $c\in Q_3(p)$, and $-\f 1c\in Q_3(p)$ if and only if $c\in
Q_1(p)$. Thus, replacing $c$ with $-1/c$ in Theorem 3.1 we deduce
the result.
 \pro{Corollary 3.2} Let $p$ be an odd prime, $c\in\Bbb
Z_p$ and $\sls {c^2+1}p=-1$. Then
$$\align 4\sum_{k=1}^{[\f{p+1}4]}\b{4k-2}{2k-1} \f 1{(16(c^2+1))^k}
&\e -2\sum_{k=(p+1)/2}^{[3p/4]}\b{4k}{2k}\f 1{(16(c^2+1))^k}
\\&\e
\sum_{k=0}^{[p/4]}\b{4k}{2k}\f 1{(16(c^2+1))^k} \mod
p.\endalign$$\endpro Proof. This is immediate from Theorem 3.1 and
(2.4).
 \pro{Corollary 3.3} Let $p$ be an odd prime, $c\in\Bbb Z_p$
and $c(c^2+1)\not\e 0\mod p$. Then
$$\sum_{k=(p+1)/2}^{[3p/4]}\b{4k}{2k}\Ls {c^2}
{16(c^2+1)}^k\e -c^2\Ls 2p\sum_{k=(p+1)/2}^{[3p/4]}\b{4k}{2k}\f
1{(16(c^2+1))^k}\mod p$$ and
$$\sum_{k=1}^{[\f{p+1}4]}\b{4k-2}{2k-1}
\Ls {c^2}{16(c^2+1)}^k\e -c^2\Ls 2p
\sum_{k=1}^{[\f{p+1}4]}\b{4k-2}{2k-1} \f 1{(16(c^2+1))^k}\mod p.$$
\endpro
Proof. This is immediate from Theorem 3.1 and Corollary 3.1.
\pro{Corollary 3.4} Let $p$ be an odd prime. Then
$$\align 2\sum_{k=(p+1)/2}^{[3p/4]}\b{4k}{2k}\f 1{32^k}
&\e -4\sum_{k=1}^{[\f{p+1}4]}\b{4k-2}{2k-1}\f 1{32^k}
\\&\e\cases 0\mod p&\t{if $p\e \pm 1\mod 8$,}
\\1\mod p&\t{if $p\e \pm 5\mod{16}$,}
\\-1\mod p&\t{if $p\e \pm 3\mod{16}$.}
\endcases\endalign$$
\endpro
Proof. Taking $c=1$ in Theorem 3.1 and then applying (2.7) we obtain
the result.
 \pro{Corollary 3.5} Let $q$ be a prime of the form $4m+1$
and $q=a^2+b^2$ with $a,b\in\Bbb Z$ and $a\e 1\mod 2$. Suppose that
$p$ is an odd prime and $\sls pq=-1$. Then
$$(-1)^{\f{p-1}2\cdot\f{q-1}4}\sum_{k=1}^{[\f{p+1}4]}
\b{4k-2}{2k-1}\Ls{b^2} {16q}^k\e \mp \f b{4a}\mod p\iff
p^{\f{q-1}4}\e\pm \f ba\mod q.$$
\endpro
Proof. As $q^{\f{p-1}2}\e \sls qp=\sls pq=-1\mod p$ we see that
$p\nmid ab$.  Taking $c=\f ab$ in Corollary 3.2 we see that
$$4\sum_{k=1}^{[\f{p+1}4]}\b{4k-2}{2k-1}\Ls{b^2}{16q}^k
\e \sum_{k=0}^{[p/4]}\b{4k}{2k}\Ls{b^2}{16q}^k\mod p.\tag 3.1$$
 Now
applying Theorem 2.3 we obtain the result.

\pro{Corollary 3.6} Let $q$ be a prime of the form $4m+1$ and
$q=a^2+b^2$ with $a,b\in\Bbb Z$ and $a\e 1\mod 2$. Suppose that $p$
is an odd prime and $\sls pq=-1$. Then
$$(-1)^{\f{p^2-1}8+\f{p-1}2\cdot\f{q-1}4}\sum_{k=1}^{[\f{p+1}4]}
\b{4k-2}{2k-1}\Ls{a^2} {16q}^k\e \mp \f a{4b}\mod p\iff
p^{\f{q-1}4}\e\pm \f ab\mod q.$$
\endpro
Proof. As $q^{\f{p-1}2}\e \sls qp=\sls pq=-1\mod p$ we see that
$p\nmid ab$.  Taking $c=\f ba$ in Corollary 3.2 we see that
$$4\sum_{k=1}^{[\f{p+1}4]}\b{4k-2}{2k-1}\Ls{a^2}{16q}^k
\e \sum_{k=0}^{[p/4]}\b{4k}{2k}\Ls{a^2}{16q}^k\mod p.\tag 3.2$$
 Now
applying Theorem 2.2 we obtain the result.
  \pro{Corollary 3.7} Let $p$ be an odd prime with $p\not=5$. Then
$$\sum_{k=1}^{[\f{p+1}4]}\b{4k-2}{2k-1}\f 1{20^k}
\e \cases 0\mod p&\t{if $p\e \pm 1\mod 5$,}
\\\mp (-1)^{\f{p-1}2}\f
12\mod p&\t{if $p\e \pm 2\mod 5$.}\endcases$$
\endpro
Proof. When $p\e \pm 1\mod 5$, taking $c=\f 12$ in Theorem 3.1 we
obtain the result. When $p\e \pm 2\mod 5$, taking $q=5$, $a=1$ and
$b=2$ in Corollary 3.5 we obtain the result.

\end{document}